\documentclass[10pt]{amsart}
\usepackage{amsfonts}
\usepackage{amsmath}
\usepackage{amsthm}
\usepackage{amssymb}
\usepackage{latexsym}

\newtheorem{lem}{Lemma}[section]

\newtheorem{prop}{Proposition}[section]

\theoremstyle{definition}

\theoremstyle{definition}
\newtheorem{thm}{Theorem}

\newtheorem*{rem}{Remark}

\newenvironment{pf}{\proof}{\endproof}

\theoremstyle{remark}

\numberwithin{equation}{section}
\begin{document}

\newcommand{\thmref}[1]{Theorem~\ref{#1}}
\newcommand{\secref}[1]{Sect.~\ref{#1}}
\newcommand{\lemref}[1]{Lemma~\ref{#1}}
\newcommand{\propref}[1]{Proposition~\ref{#1}}
\newcommand{\corref}[1]{Corollary~\ref{#1}}
\newcommand{\remref}[1]{Remark~\ref{#1}}
\newcommand{\er}[1]{(\ref{#1})}
\newcommand{\nc}{\newcommand}
\newcommand{\rnc}{\renewcommand}
\nc{\goth}{\mathfrak}
\rnc{\bold}{\mathbf}
\renewcommand{\frak}{\mathfrak}
\renewcommand{\Bbb}{\mathbb}

\nc{\Cal}{\mathcal}
\nc{\Xp}[1]{X^+(#1)}
\nc{\Xm}[1]{X^-(#1)}
\nc{\on}{\operatorname}
\nc{\ch}{\mbox{ch}}
\nc{\Z}{{\bold Z}}
\nc{\J}{{\mathcal J}}
\nc{\C}{{\bold C}}
\nc{\Q}{{\bold Q}}
\renewcommand{\P}{{\mathcal P}}
\nc{\N}{{\Bbb N}}
\nc\beq{\begin{equation}}
\nc\enq{\end{equation}}
\nc\lan{\langle}
\nc\ran{\rangle}
\nc\bsl{\backslash}
\nc\mto{\mapsto}
\nc\lra{\leftrightarrow}
\nc\hra{\hookrightarrow}
\nc\sm{\smallmatrix}
\nc\esm{\endsmallmatrix}
\nc\sub{\subset}
\nc\ti{\tilde}
\nc\nl{\newline}
\nc\fra{\frac}
\nc\und{\underline}
\nc\ov{\overline}
\nc\ot{\otimes}
\nc\bbq{\bar{\bq}_l}
\nc\bcc{\thickfracwithdelims[]\thickness0}
\nc\ad{\text{\rm ad}}
\nc\Ad{\text{\rm Ad}}
\nc\Hom{\text{\rm Hom}}
\nc\End{\text{\rm End}}
\nc\Ind{\text{\rm Ind}}
\nc\Res{\text{\rm Res}}
\nc\Ker{\text{\rm Ker}}
\rnc\Im{\text{Im}}
\nc\sgn{\text{\rm sgn}}
\nc\tr{\text{\rm tr}}
\nc\Tr{\text{\rm Tr}}
\nc\supp{\text{\rm supp}}
\nc\card{\text{\rm card}}
\nc\bst{{}^\bigstar\!}
\nc\he{\heartsuit}
\nc\clu{\clubsuit}
\nc\spa{\spadesuit}
\nc\di{\diamond}

\nc\al{\alpha}
\nc\bet{\beta}
\nc\ga{\gamma}
\nc\de{\delta}
\nc\ep{\epsilon}
\nc\io{\iota}
\nc\om{\omega}
\nc\si{\sigma}
\rnc\th{\theta}
\nc\ka{\kappa}
\nc\la{\lambda}
\nc\ze{\zeta}

\nc\vp{\varpi}
\nc\vt{\vartheta}
\nc\vr{\varrho}

\nc\Ga{\Gamma}
\nc\De{\Delta}
\nc\Om{\Omega}
\nc\Si{\Sigma}
\nc\Th{\Theta}
\nc\La{\Lambda}

\nc\boa{\bold a}
\nc\bob{\bold b}
\nc\boc{\bold c}
\nc\bod{\bold d}
\nc\boe{\bold e}
\nc\bof{\bold f}
\nc\bog{\bold g}
\nc\boh{\bold h}
\nc\boi{\bold i}
\nc\boj{\bold j}
\nc\bok{\bold k}
\nc\bol{\bold l}
\nc\bom{\bold m}
\nc\bon{\bold n}
\nc\boo{\bold o}
\nc\bop{\bold p}
\nc\boq{\bold q}
\nc\bor{\bold r}
\nc\bos{\bold s}
\nc\bou{\bold u}
\nc\bov{\bold v}
\nc\bow{\bold w}
\nc\boz{\bold z}

\nc\ba{\bold A}
\nc\bb{\bold B}
\nc\bc{\bold C}
\nc\bd{\bold D}
\nc\be{\bold E}
\nc\bg{\bold G}
\nc\bh{\bold H}
\nc\bi{\bold I}
\nc\bj{\bold J}
\nc\bk{\bold K}
\nc\bl{\bold L}
\nc\bm{\bold M}
\nc\bn{\bold N}
\nc\bo{\bold O}
\nc\bp{\bold P}
\nc\bq{\bold Q}
\nc\br{\bold R}
\nc\bs{\bold S}
\nc\bt{\bold T}
\nc\bu{\bold U}
\nc\bv{\bold V}
\nc\bw{\bold W}
\nc\bz{\bold Z}
\nc\bx{\bold X}

\nc\ca{\mathcal A}
\nc\cb{\mathcal B}
\nc\cc{\mathcal C}
\nc\cd{\mathcal D}
\nc\ce{\mathcal E}
\nc\cf{\mathcal F}
\nc\cg{\mathcal G}
\rnc\ch{\mathcal H}
\nc\ci{\mathcal I}
\nc\cj{\mathcal J}
\nc\ck{\mathcal K}
\nc\cl{\mathcal L}
\nc\cm{\mathcal M}
\nc\cn{\mathcal N}
\nc\co{\mathcal O}
\nc\cp{\mathcal P}
\nc\cq{\mathcal Q}
\nc\car{\mathcal R}
\nc\cs{\mathcal S}
\nc\ct{\mathcal T}
\nc\cu{\mathcal U}
\nc\cv{\mathcal V}
\nc\cz{\mathcal Z}
\nc\cx{\mathcal X}
\nc\cy{\mathcal Y}

\nc\e[1]{E_{#1}}
\nc\ei[1]{E_{\delta - \alpha_{#1}}}
\nc\esi[1]{E_{s \delta - \alpha_{#1}}}
\nc\eri[1]{E_{r \delta - \alpha_{#1}}}
\nc\ed[2][]{E_{#1 \delta,#2}}
\nc\ekd[1]{E_{k \delta,#1}}
\nc\emd[1]{E_{m \delta,#1}}
\nc\erd[1]{E_{r \delta,#1}}

\nc\ef[1]{F_{#1}}
\nc\efi[1]{F_{\delta - \alpha_{#1}}}
\nc\efsi[1]{F_{s \delta - \alpha_{#1}}}
\nc\efri[1]{F_{r \delta - \alpha_{#1}}}
\nc\efd[2][]{F_{#1 \delta,#2}}
\nc\efkd[1]{F_{k \delta,#1}}
\nc\efmd[1]{F_{m \delta,#1}}
\nc\efrd[1]{F_{r \delta,#1}}

\nc\fa{\frak a}
\nc\fb{\frak b}
\nc\fc{\frak c}
\nc\fd{\frak d}
\nc\fe{\frak e}
\nc\ff{\frak f}
\nc\fg{\frak g}
\nc\fh{\frak h}
\nc\fj{\frak j}
\nc\fk{\frak k}
\nc\fl{\frak l}
\nc\fm{\frak m}
\nc\fn{\frak n}
\nc\fo{\frak o}
\nc\fp{\frak p}
\nc\fq{\frak q}
\nc\fr{\frak r}
\nc\fs{\frak s}
\nc\ft{\frak t}
\nc\fu{\frak u}
\nc\fv{\frak v}
\nc\fz{\frak z}
\nc\fx{\frak x}
\nc\fy{\frak y}

\nc\fA{\frak A}
\nc\fB{\frak B}
\nc\fC{\frak C}
\nc\fD{\frak D}
\nc\fE{\frak E}
\nc\fF{\frak F}
\nc\fG{\frak G}
\nc\fH{\frak H}
\nc\fJ{\frak J}
\nc\fK{\frak K}
\nc\fL{\frak L}
\nc\fM{\frak M}
\nc\fN{\frak N}
\nc\fO{\frak O}
\nc\fP{\frak P}
\nc\fQ{\frak Q}
\nc\fR{\frak R}
\nc\fS{\frak S}
\nc\fT{\frak T}
\nc\fU{\frak U}
\nc\fV{\frak V}
\nc\fZ{\frak Z}
\nc\fX{\frak X}
\nc\fY{\frak Y}
\nc\tfi{\ti{\Phi}}
\nc\bF{\bold F}
\rnc\bol{\bold 1}

\nc\ua{\bold U_\A}

\nc\qinti[1]{[#1]_i}
\nc\q[1]{[#1]_q}
\nc\xpm[2]{E_{#2 \delta \pm \alpha_#1}}  
\nc\xmp[2]{E_{#2 \delta \mp \alpha_#1}}
\nc\xp[2]{E_{#2 \delta + \alpha_{#1}}}
\nc\xm[2]{E_{#2 \delta - \alpha_{#1}}}
\nc\hik{\ed{k}{i}}
\nc\hjl{\ed{l}{j}}
\nc\qcoeff[3]{\left[ \begin{smallmatrix} {#1}& \\ {#2}& \end{smallmatrix}
\negthickspace \right]_{#3}}
\nc\qi{q}
\nc\qj{q}

\nc\ufdm{{_\ca\bu}_{\rm fd}^{\le 0}}


\nc\isom{\cong} 

\nc{\pone}{{\Bbb C}{\Bbb P}^1}
\nc{\pa}{\partial}
\def\H{\mathcal H}
\def\L{\mathcal L}
\nc{\F}{{\mathcal F}}
\nc{\Sym}{{\goth S}}
\nc{\A}{{\mathcal A}}
\nc{\arr}{\rightarrow}
\nc{\larr}{\longrightarrow}

\nc{\ri}{\rangle}
\nc{\lef}{\langle}
\nc{\W}{{\mathcal W}}
\nc{\uqatwoatone}{{U_{q,1}}(\su)}
\nc{\uqtwo}{U_q(\goth{sl}_2)}
\nc{\dij}{\delta_{ij}}
\nc{\divei}{E_{\alpha_i}^{(n)}}
\nc{\divfi}{F_{\alpha_i}^{(n)}}
\nc{\Lzero}{\Lambda_0}
\nc{\Lone}{\Lambda_1}
\nc{\ve}{\varepsilon}
\nc{\phioneminusi}{\Phi^{(1-i,i)}}
\nc{\phioneminusistar}{\Phi^{* (1-i,i)}}
\nc{\phii}{\Phi^{(i,1-i)}}
\nc{\Li}{\Lambda_i}
\nc{\Loneminusi}{\Lambda_{1-i}}
\nc{\vtimesz}{v_\ve \otimes z^m}

\nc{\asltwo}{\widehat{\goth{sl}_2}}
\nc\ag{\widehat{\goth{g}}}  
\nc\teb{\tilde E_\boc}
\nc\tebp{\tilde E_{\boc'}}

\title[Quantum affine algebras at roots of unity]
{Realization of level one representations of
$U_q(\widehat{\mathfrak g})$ at a root of unity}
\author{Vyjayanthi Chari}
\address{Vyjayanthi Chari, University of California, Riverside}
\email{chari@math.ucr.edu}
\author{Naihuan Jing}
\address{Naihuan Jing, North Carolina State University, Raleigh}
\email{jing@math.ncsu.edu}
\thanks{N.J. is partially supported by NSF grant
DMS-9970493}

\begin{abstract}
Using vertex operators, we construct explicitly Lusztig's $\mathbb Z[q, q^{-1}]$-lattice for
the level one irreducible representations of quantum affine algebras of ADE type. We then realize the level one irreducible modules at roots of unity and show that
the $q$-dimension is still given by the Weyl-Kac character formula. As a consequence we also answer the corresponding
question of realizing the affine Kac-Moody Lie algebras
of simply laced type at level one in finite characteristic.
\end{abstract}

\maketitle

\section{Introduction}
In \cite{L, L1} Lusztig proved that a quantum Kac--Moody algebra $\bu$ defined over $\bq(q)$ admits an
$\A=\bz[q,q^{-1}]$--lattice $\ua$ and that any irreducible 
highest weight integrable representation $V$ of $\bu$ admits a  corresponding $\A$--lattice, say $V_\A$. This allows us to specialize  $q$ to a non--zero complex number $\zeta$ and we let $\bu_\zeta$, $W_\zeta$ denote the corresponding objects.
If $\zeta$ is not a root of unity, Lusztig proved that $W_\zeta$ is irreducible and its character is the same as that of the corresponding classical representation. On the other hand, when $\zeta$ is a primitive $l^{th}$ root of unity, the situation is more interesting, even for finite--dimensional Kac--Moody algebras. In that case, $W_\zeta$ is not always irreducible: a sufficient condition for irreducibility \cite{APW} is that the highest weight $\Lambda$ of $V$ should be \lq\lq small\rq\rq \ in the sense
that $(\Lambda,\alpha)<l$ for all positive roots $\alpha$. The corresponding question for infinite--dimensional Kac--Moody algebras at roots of unity is open, and in  this paper we  answer it  in the case of level one representations of quantum affine algebras of ADE type. Note that the condition $(\Lambda,\alpha)<l$ never holds in this case; nevertheless, we find that $W_\zeta$ is irreducible provided that $l$ is coprime to the Coxeter number of the underlying 
finite--dimensional Lie algebra.

 The level one  representations of an affine Lie algebra of ADE type can be explicitly constructed in the tensor product of  a symmetric algebra and a twisted group algebra \cite{FK, S}. Essentially, these 
representations are built from the canonical representation of an infinite--dimensional Heisenberg algebra.  Later, in \cite{FJ} this construction was extended to the case of the basic representations of the quantum affine algebras of ADE type. Again, the representations are built from the representation of a 
suitable quantum Heisenberg algebra.  In this paper, we identify the natural lattice $V_\A$ of the level one representation explicitly as the tenosr product of the lattice of Schur functions tensored with the obvious $\A$--lattice in the twisted group algebra (see also \cite{J3}).  We also  describe the action of the divided powers of the 
Chevalley (and Drinfeld) generators on an $\A$--basis of $V_\A$ and  this  allows us to realize the level one irreducible
representations $W_\zeta$ explicitly  and prove that they are irreducible.

Our methods also apply to the study of highest weight representations of affine Lie algebras in characteristic $p$, and the corresponding results are also new in that situation. In particular we give
an explicit realization of the $\mathbb Z$-form \cite{Br} of the vertex
representation of the affine Lie algebras.

\section{ The algebras $\bu$, $\ua$.} Throughout this paper $\frak{g}$
will denote a simply-laced, finite-dimensional complex simple Lie algebra
and $(a_{ij})_{i,j\in I}, \ I=\{1,\dots, n\},$ will denote its Cartan
matrix.  Let $(a_{ij})_{i,j\in\hat I}, \ \hat I=I \cup \{0\}$, be the
extended Cartan matrix of $\frak g$ and let $\ag$ be the corresponding
affine Lie algebra. Let $R$ (resp. $R^+$) denote a set of roots (resp. 
positive roots) of $\frak g$ and let $\alpha_i$ ($i\in I$) be a set of
simple roots. Let $Q$ be the root lattice of $\frak g$, $P$ the weight
lattice and let $\omega_i\in P$ ($i\in I$) be the fundamental weights of
$\frak g$.  For $\omega\in P$, $\eta\in Q$, define an integer
$|\omega|\cdot|\eta|$ by extending bilinearly the assignment
$|\omega_i|\cdot|\alpha_j| =\delta_{ij}$.  Notice that
$|\alpha_i|\cdot|\alpha_j| =a_{ij}$. Let $\theta$ be the highest root of
$\frak g$. 

Let $q$ be an indeterminate, let $\mathbb{Q}(q)$ be the field of rational
functions in $q$ with rational coefficients, and let $\mathcal{A}=\mathbb
{Z}[q,q^{-1}]$ be
the ring of Laurent polynomials with integer coefficients. For
$r,m\in\mathbb N$, $m\ge r$, define
\begin{equation*}
[m]=\frac{q^m -q^{-m}}{q -q^{-1}},\ \ \ \
  [m]! =[m][m-1]\ldots [2][1],\ \ \ \
\left[\begin{matrix} m\\
  r\end{matrix}\right] = \frac{[m]!}{[r]![m-r]!}.
\end{equation*}
  Then $\left[\begin{matrix} m\\r\end{matrix}\right]\in\mathcal A$
  for all $m\ge r\ge 0$.

\begin{prop}{\label{defnbu}} There is a Hopf algebra $\bu$ over $\bq(q)$
which is generated as an algebra by elements $E_{\alpha_i}$,
$F_{\alpha_i}$, $K_i^{{}\pm 1}$ ($i\in\hat I$), $D^{\pm 1}$ with the
following defining relations: \begin{align*}
  K_iK_i^{-1}=K_i^{-1}K_i&=1,\ \ \ \ K_iK_j=K_jK_i,\\
  K_iD=DK_i,\ \  DD^{-1}&=D^{-1}D=1,\\
  DE_{\alpha_i}D^{-1}=q^{\delta_{i0}}E_{\alpha_i},\ \ 
  &DF_{\alpha_i}D^{-1}=q^{-\delta_{i0}}F_{\alpha_i}, \\ 
  K_iE_{\alpha_j} K_i^{-1}&=q^{ a_{ij}}E_{\alpha_j},\\ 
K_iF_{\alpha_j} K_i^{-1}&=q^{-a_{ij}}F_{\alpha_j},\\
  [E_{\alpha_i}, F_{\alpha_j}
]&=\delta_{ij}\frac{K_i-K_i^{-1}}{q-q^{-1}},\\ 
  \sum_{r=0}^{1-a_{ij}}(-1)^r\left[\begin{matrix} 1-a_{ij}\\ 
  r\end{matrix}\right]
&(E_{\alpha_i})^rE_{\alpha_j}(E_{\alpha_i})^{1-a_{ij}-r}=0\ 
  \ \ \ \ \text{if $i\ne j$},\\
\sum_{r=0}^{1-a_{ij}}(-1)^r\left[\begin{matrix} 1-a_{ij}\\ 
  r\end{matrix}\right]
&(F_{\alpha_i})^rF_{\alpha_j}(F_{\alpha_i})^{1-a_{ij}-r}=0\ 
  \ \ \ \ \text{if $i\ne j$}. \end{align*} The comultiplication of $\bu$
is given on generators by \begin{align*}
\Delta(E_{\alpha_i})&=E_{\alpha_i}\ot 1+K_i\ot E_{\alpha_i},\ \
\Delta(F_{\alpha_i})=F_{\alpha_i}\ot K_i^{-1} + 1\ot F_{\alpha_i},\\
\Delta(K_i)&=K_i\ot K_i, \qquad \Delta(D)=D\ot D, \end{align*} for
$i\in\hat I$.\hfill\qedsymbol \end{prop} Let $\bu^+$ (resp. $\bu^-$;
$\bu^0$) be the $\bq(q)$-subalgebras of $\bu$ generated by the
$E_{\alpha_i}$ (resp. $F_{\alpha_i}$; $K_i^{\pm 1}$ and $D^{\pm 1}$) for
$i\in \hat{I}$. The following result is well-known, see \cite{L} for
instance. 

\begin{lem}{\label {butriangle}}
 $\bu\isom \bu^-\otimes\bu^0\otimes\bu^+$ as $\bq(q)$-vector
spaces.\hfill\qedsymbol \end{lem}

 It is convenient to use the following notation:
\begin{equation*}E_{\alpha_i}^{(r)}=\frac{E_{\alpha_i}^r}{[r]!}.\end{equation*}
The elements $F_{\alpha_i}^{(r)}$ are defined similarly. Let $\ua$ denote
the $\A$--subalgebra of $\bu$ generated by $E_{\alpha_i}^{(r)}$,
$F_{\alpha_i}^{(r)}$, $K_i^{\pm 1}$ ($i\in \hat I$) and $D^{\pm 1}$. The
subalgebras $\ua^\pm$ are defined in the obvious way. 

For $i\in \hat{I}$, $r\ge 1$, $m\in\mathbb Z$, define elements
\begin{align*}\genfrac{[}{]}{0pt}{}{K_i,m}{r}& =
 \prod_{s=1}^r \frac{K_i q^{m-s+1} - K_i^{-1}
    q^{-m+ s-1}}{q^s - q^{-s}},\\
\genfrac{[}{]}{0pt}{}{D,m}{r}& =
 \prod_{s=1}^r \frac{D q^{m-s+1} - D^{-1}
    q^{-m+ s-1}}{q^s - q^{-s}}.\\
\end{align*}

Let $\ua^0$ be the $\mathcal A$--subalgebra of $\ua$ generated by
$K_i^{\pm 1}$, $D^{\pm 1}$, $\genfrac{[}{]}{0pt}{}{K_i,m}{r}$ and
$\genfrac{[}{]}{0pt}{}{D,m}{r}$, $i\in\hat{I}$, $r\ge 1$ and $m\in\mathbb
Z$.  The following is well--known (see \cite{L}).
\begin{lem}{\label{atriangle}} We have $\ua\cong\ua^-\ua^0\ua^+$.\end{lem}

We shall also need another realization of $\bu$, due to \cite{Dr, B, J2}.
\begin{thm}{\label{newr}} There is an isomorphism 
  of $\bq(q)$-Hopf algebras from $\bu$ to the algebra with generators
  $x_{i,r}^{{}\pm{}}$ ($i\in I$, $r\in\bz$), $K_i^{{}\pm 1}$ ($i\in I$),
  $h_{i,r}$ ($i\in I$, $r\in \bz\backslash\{0\}$) and $C^{{}\pm 1}$,
  and the following defining relations:
\begin{align*}
  C^{\pm 1}\ &\text{are central,}\\ K_iK_i^{-1} = K_i^{-1}K_i
  =1,\;\; &CC^{-1} =C^{-1}C =1,\\ K_iK_j =K_jK_i,\;\;
  &K_ih_{j,r} =h_{j,r}K_i,\\  K_ix_{j,r}^\pm K_i^{-1} &= q^{{}\pm
    a_{ij}}x_{j,r}^{{}\pm{}},\\
DD^{-1} =D^{-1}D =1,\;\;& DK_i =K_iD,\\
Dh_{j,r}D^{-1} =q^rh_{j,r},\;\;& Dx_{j,r}^{\pm 1}D^{-1} = q^rx_{j,r}^\pm,\\ 
  [h_{i,r},h_{j,s}]&=\delta_{r,-s}\frac1{r}[ra_{ij}]\frac{C^r-C^{-r}}
  {q-q^{-1}},\\ 
[h_{i,\pm r} , x_{j,s}^{{}\pm{}}] &=
  \pm\frac1r[ra_{ij}]x_{j,s\pm r}^{{}\pm{}},\ \ r>0, \\ 
[h_{i,\mp r} , x_{j,s}^{{}\pm{}}] &=
  \pm\frac1rC^r[ra_{ij}]x_{j,s\pm r}^{{}\pm{}},\ \ r>0, \\ 
x_{i,r+1}^{{}\pm{}}x_{j,s}^{{}\pm{}} -q^{{}\pm
    a_{ij}}x_{j,s}^{{}\pm{}}x_{i,r+1}^{{}\pm{}} &=q^{{}\pm
    a_{ij}}x_{i,r}^{{}\pm{}}x_{j,s+1}^{{}\pm{}}
  -x_{j,s+1}^{{}\pm{}}x_{i,r}^{{}\pm{}},\\ [x_{i,r}^+ ,
  x_{j,s}^-]=\delta_{i,j} & \frac{ C^{-s}\psi_{i,r+s}^+ -
    C^{-r} \psi_{i,r+s}^-}{q - q^{-1}},\\ 
\sum_{\pi\in\Sigma_m}\sum_{k=0}^m(-1)^k\left[\begin{matrix}m\\k\end{matrix}
\right]
  x_{i, r_{\pi(1)}}^{{}\pm{}}\ldots x_{i,r_{\pi(k)}}^{{}\pm{}} &
  x_{j,s}^{{}\pm{}} x_{i, r_{\pi(k+1)}}^{{}\pm{}}\ldots
  x_{i,r_{\pi(m)}}^{{}\pm{}} =0,\ \ \text{if $i\ne j$},
\end{align*}
for all sequences of integers $r_1,\ldots, r_m$, where $m =1-a_{ij}$, $\Sigma_m$ is the symmetric group on $m$ letters, and the $\psi_{i,r}^{{}\pm{}}$ are determined by equating powers of $u$ in the formal power series 
$$\sum_{r=0}^{\infty}\psi_{i,\pm r}^{{}\pm{}}u^{{}\pm r} = K_i^{{}\pm 1} exp\left(\pm(q-q^{-1})\sum_{s=1}^{\infty}h_{i,\pm s} u^{{}\pm s}\right).$$

\hfill\qedsymbol\end{thm}

Following \cite[Section 3]{CP},
we define elements  $P_{k,i}$ and $\tilde{P}_{k,i}$ via  the generating functions
\begin{equation} \label{integralimaginary}
 \mathcal{P}_i^\pm(u)= \sum_{k \ge 0}{ P^\pm_{i,k}} u^k = \exp \left(- \sum_{k = 1}^\infty
  \frac{h_{i,\pm k}}{[k]}u^k\right)= \exp \left( -\sum_{k = 1}^\infty
  \frac{\tilde{h}_{i,\pm k}}{k}u^k\right), 
\end{equation} 
\begin{equation} \label{tildeP} \tilde{\mathcal{P}}_i^\pm(u)=\sum_{k \ge 0}{\tilde{P}^\pm_{i, k}} u^k = 
  \exp \left( \sum_{k = 1}^\infty \frac{h_{i,\pm k}}{[k]}u^k\right)=\exp \left( \sum_{k = 1}^\infty
  \frac{\tilde{h}_{i,\pm k} }{k}u^k\right),
\end{equation}
where $\tilde{h}_{i,k} = \frac{kh_{i,k}}{[k]}$. Notice that these formulas  are exactly those that relate  the elementary symmetric functions (resp. complete symmetric functions) to the power sum symmetric functions \cite{M}.  For a vertex operator approach to this and to Schur functions, see \cite{J3}.

The following result was proved in \cite[Section 5]{CP}.
\begin{lem}  For all $i\in I$, $k\in\bz$, $k\ge 0$, we have
\begin{equation*} P^\pm_{i,k},\ \ \tilde{P}^\pm_{i, k}\in\ua.\end{equation*}\end{lem}

Let $\tilde\ua$ (resp. $\tilde\ua^\pm$) be the $\mathcal A$--subalgebra generated by $(x_{i,n}^\pm )^{(r)}$, $r,n\in\bz$, $r\ge 0$, $i\in I$, (resp. $n\in\bz$, $\pm n\ge 0$) and $\ua^0$.
 The following result is proved in \cite[Section 2]{BCP}.
\begin{prop}{\label{adrin}} We have,
\begin{equation*}\ua=\tilde\ua,\ \ \ \ \ua^\pm\subset\tilde\ua^\pm. 
\ \ \ \ \ \ \qedsymbol\end{equation*}

\end{prop}

Finally, let $\bu(0)$ (resp. $\ua(0)$) be the $\bq(q)$--subalgebra of $\bu$ (resp. the $\A$--subalgebra of $\ua$) generated by the elements $h_{i,n}$, $i\in I$, $n\in\bz$ (resp. $P^\pm_{i,k}$, $i\in I$, $k\in\bz$, $k>0$), $C^{\pm 1}$. The 
subalgebras $
\bu^\pm(0)$ and $\ua^\pm(0)$ are defined in the obvious way.
\begin{prop} \label{heis}

\begin{enumerate} 

\item[(i)] The algebra $\bu(0)$ is defined by the relations   
\begin{align*} [h_{i,n}, h_{j,m}]& = \frac 1n\delta_{m,-n}{[na_{ij}]}\frac{C^n-C^{-n}}{q-q^{-1}},\\
C^{\pm 1}h_{i,n}& = h_{i,n}C^{\pm 1},\\
\end{align*}
for all $i,j\in I$ and $m,n\in\bz$. In particular, $\bu^\pm(0)$ is commutative.
\item[(ii)] For all $i\in I$, $k>0$, we have
\begin{align*}
P^\pm _{i,k}&=-\frac 1k\sum_{m=0}^k \tilde{h}_{i,m}P^\pm_{i,k-m},\\
\tilde{P}^\pm _{i,k}&=\frac 1k\sum_{m=0}^k \tilde{h}_{i,m}\tilde{P}^\pm_{i,k-m}.\end{align*}
In particular, $\tilde{h}_{i,k}, \tilde{P}^\pm_{i,k}\in\ua(0)$ and 
as $\bq(q)$--spaces we have
\begin{align*}
\bu(0)\cong \bq(q)\otimes_\A\ua(0),\\
\bu^\pm(0)\cong \bq(q)\otimes_\A\ua^\pm(0).\end{align*}\end{enumerate}
\item[(iii)]Monomials in $P^\pm_{i,  n}$ (resp. $\tilde{P}^\pm_{i, n}$), $i\in I$, $n>0$, form a basis for $\ua^\pm(0)$. \end{prop}
\begin{pf} Part (i) is a consequence of the PBW theorem for
$\bu$ proved in \cite{B}. Parts (ii) and (iii) follow from the definition of the elements $\tilde{P}^\pm_{i,k}$ (see \cite{BCP} for details).
\end{pf}

\section{ The level one representations of $\bu$ and $\ua$} 
We begin this section by recalling  the natural irreducible representation of $\bu(0)$ and 
we construct a natural $\ua(0)$--lattice in this representation. We then recall the definition of the highest weight representations $V_q(\Lambda)$ of $\bu$ and   the lattice $V_{\mathcal A}(\Lambda)$ of $\ua$, see \cite{L}. Finally, we recall  the explicit construction of the level one representations given in \cite{FJ} and state and prove the main theorem of the paper.

Consider the left ideal $\mathcal I$ in $\bu(0)$ generated by $C^{\pm 1}-q^{\pm 1}$ and $\bu^+(0)$. Then, $\bu(0)/\mathcal I$ is a left $\bu(0)$--module through left multiplication.  It is  easy to see that as $\bq(q)$--spaces, we have
$$ \bu^-(0)\cong \bu(0) /\mathcal I.$$
Thus $ \bu^-(0)$ acquires the structure of a left $\bu(0)$--module, and 
we 
let $\pi:\bu(0)\to{\text{End}}(\bu^-(0))$ be this representation. Then, elements of $\bu^-(0)$ act by left multiplication and it is easy to see that for $n>0$, $i\in I$, $\pi(h_{i,n})$ is the derivation of $\bu^-(0)$ obtained by extending the assignment,
\begin{equation*}\pi(h_{i,n})h_{j,-m} = \delta_{n,m}\frac{[na_{ij}][n]}{n}.\end{equation*}

\begin{prop} {\label{heisr}}

\begin{enumerate} 

\item[(i)] $\pi$ is an irreducible representation of $\bu(0)$. 
\item[(ii)]For $i,j\in I$, we have,
\begin{equation*}\pi(\tilde\mathcal{P}_i^+(u)).\tilde\mathcal{P}_j^-(v) = f_{i,j}(u,v)\tilde\mathcal{P}_j^-(v),\end{equation*}
\begin{equation*}\pi(\mathcal{P}_i^+(u)).\mathcal{P}_j^-(v)= f_{i,j}(u,v)\mathcal{P}_j^-(v),\end{equation*}
\begin{equation*}f_{i,j}(u,v)\pi(\mathcal{P}_i^+(u)).\tilde\mathcal{P}_j^-(v)= \tilde\mathcal{P}_j^-(v),\end{equation*}
where the power series $f_{i, j}$ is defined by
\begin{align*} {f}_{i,j}(u,v)&= 1\ \ \ \text{if $a_{ij}=0$},\\
&=(1-uv)\ \ \  \text{if $a_{ij}=-1$},\\
&=(1-quv)^{-1}(1-q^{-1}uv)^{-1}\ \ \  \text{if $a_{ij}=2$}.\end{align*}
\item[(iii)] $\pi(\ua(0)) \ua^-(0)\subset \ua^-(0)$.
\end{enumerate}
\end{prop}
\begin{proof} Part (i) is well-known.
For (ii), notice that the relations in Proposition \ref{heis} imply that
\begin{align*}&\pi(\tilde\mathcal{P}_-^+(u)).\tilde\mathcal{P}_j^-(v)\\
&=\exp \left( \sum_{k = 1}^\infty
  \frac{\pi(\tilde{h}_{i, k})} {k}u^k\right) \exp \left( \sum_{k = 1}^\infty
  \frac{\tilde{h}_{j,-k}}{k}u^k\right)\\
&=\exp\left(\sum_{k=1}^\infty\frac{[ka_{ij}]}{k[k]}u^kv^k\right)\exp \left( \sum_{k = 1}^\infty
  \frac{\tilde{h}_{j,-k}}{k}u^k\right)\exp\left(\sum_{k = 1}^\infty
  \frac{\pi(\tilde{h}_{i, k})} {k}u^k \right).1\\
&=f_{i,j}(u,v)\tilde\mathcal{P}_j^-(v).\end{align*}
The second equality above follows by using the Campbell--Hausdorff formula.
The calculation of $f_{i,j}(u,v)$ is now straightforward. The other equations are proved similarly. 
 Part (iii)  follows immediately from (ii).
\end{proof}

By a weight, we mean a pair $(\mu, n)\in \bz^{|\hat I|}\times \bz$. If $n=0$, we shall denote the pair $(\mu,0)$ as $\mu$. 
A representation $W$ of $\bu$ is said to  be of type 1 if
$$W=\bigoplus_{(\mu, n)} W_{\mu,n},$$
where $W_{\mu,n}=\{w\in W| K_i.w=q^{\mu_i}w,\ \ D.w =q^n w\}$. If $W_{\mu,n}\ne 0$, then $W_{\mu,n}$ is called the weight space of $W$ with weight $(\mu,n)$.
Throughout this paper we will consider only type 1 representations. Writing $\theta=\sum_{i\in I}d_i\alpha_i$, we define the level of $(\mu, n)$ to be 
$\sum_{i\in I}d_i\mu_i+\mu_0$.

For $i\in\hat I$, let $\Lambda_i$ be the $\hat I$--tuple with one in the $i^{th}$ place and zero elsewhere.  Given a weight $\Lambda=\sum_in_i\Lambda_i$, $n_i\ge 0$, let $V_q(\Lambda)$ be the irreducible highest weight  $\bu$--module with highest weight
$\Lambda$ and let $v_{\Lambda}$ be the highest weight
vector. 
Thus, $V_q(\Lambda)$ is generated by $v_\Lambda$ with relations, \begin{equation*} E_{\alpha_i}.v_{\Lambda}=0,\ \
K_i.v_{\Lambda}=q^{n_i}v_{\Lambda},\ \ D.v_{\Lambda}=v_{\Lambda},
\ \ F_{\alpha_i}^{n_i+1}.v_{\Lambda}=0,\ \ \end{equation*}
for $i\in\hat{I}$. Clearly $V_q(\Lambda)$ is of type 1. We say that $V_q(\Lambda)$ has level one if $\Lambda$ has level one.

Set \begin{equation*} V_{\mathcal A}(\Lambda)=\ua.v_\Lambda.\end{equation*}By Lemma \ref{atriangle} we see that $V_\A =\ua^-.v_\Lambda$. The following result is now an immediate consequence of Proposition \ref{adrin}.
\begin{lem} We have
\begin{equation*} V_\A(\Lambda)=\tilde\ua.v_\Lambda=\tilde\ua^-.v_\Lambda.\end{equation*}
\hfill\qedsymbol\end{lem}

The following result is due to Lusztig \cite{L}.
\begin{prop}{\label{vlatt}} $V_{\mathcal A}(\Lambda)$ is a $\ua$--submodule of
$V_q(\Lambda)$ such that
\begin{equation*} V_q(\Lambda)\cong V_{\mathcal A}(\Lambda)\otimes_{\mathcal A}\bq(q).\end{equation*} Further,
\begin{equation*}V_{\mathcal A}(\Lambda)=\bigoplus_{\mu,n} V_{\mathcal A}(\Lambda)\cap
V_q(\Lambda)_{\mu,n},\end{equation*} and
\begin{equation*} {\text{dim}}_{\mathcal A}(V_{\mathcal A}(\Lambda)\cap
V_q(\Lambda)_{\mu,n}) ={\text{dim}}_{\bq(q)}
V_q(\Lambda)_{\mu,n}.\ \ \ \ \ \ \qedsymbol\end{equation*}\end{prop}

We turn now to the realization of the level one representations of $\bu$. In fact we shall restrict ourselves to constructing  the basic representation of $\hat{\frak g}$, i.e. the representation corresponding to $\Lambda_0$. The construction of the other
level one representations is  identical except
that one adjoins $v_{\Lambda_i}$ to the twisted group
algebra (see \cite{FJ}).

Fix a bilinear map $\epsilon: Q\times Q\to\{\pm 1\}$ such that for all $i\in I$, $\al, \beta, \ga\in Q$, we have,
\begin{align*}
\ep(\al, 0)&=\ep(0, \al)=1, \\
\ep(\al, \beta)\ep(\al+\beta, \ga)&=\ep(\al, \beta+\ga)\ep(\beta, \ga),\\
\epsilon(\al, \beta)\epsilon(\beta,\al)&= (-1)^{|\al|\cdot|\beta|}.
\end{align*}
Let $\bq(q)[Q]$ be the twisted group
algebra over $\bq(q)$ of the weight lattice of $\frak{g}$. Thus, $\bq(q)[Q]$ is the algebra generated by elements $e^\eta$, $\eta\in Q$, subject to  the relation, 
\begin{equation*}
e^\eta.e^{\eta'} = \epsilon(\eta,\eta')e^{\eta+\eta'}.\end{equation*}
Set
\begin{equation*}\mathcal{V}_q=\bu^-(0)\otimes \bq(q)[Q].\end{equation*}

 Let $z^\partial_i:\mathcal{V}_q\to\mathcal{V}_q[z,z^{-1}]$ be the $\bq(q)$--linear  map defined by extending \begin{equation*} z^{\partial_i}(v\otimes e^{\eta})=  (v\otimes e^{\eta})z^{|\eta|\cdot|\alpha_i|},\ \ v\in\bu^-(0),\ \ \eta\in Q.\end{equation*}

Define operators $X^\pm_{i,n}$ on $\mathcal{V}_q$ by means of the
following generating series:
\begin{align*}
X_i^{+}(z)&=\pi(\tilde{\mathcal{P}}_i^-(z))\pi(\mathcal{P}_i^+(q^{-1}z^{-1}))e^{\alpha_i}z^{\partial_i}\\
&=\sum_{n\in{\mathbb Z}}X_{i,n}^{+}z^{-n-1},\\
X_i^{-}(z)&=\pi({\mathcal{P}}_i^-(qz))\pi(\tilde{\mathcal{P}}_i^+(z^{-1}))e^{-\alpha_i}z^{-\partial_i}\\
&=\sum_{n\in\bz}X_{i,n}^{-}z^{-n-1}.
\end{align*}

The following result was proved in \cite{FJ}.
\begin{thm}{\label{realiz}} The assignment $x^\pm_{i,n}\to X^{\pm}_{i,n}$,
$h_{i,n}\to \pi(h_{i,n})\otimes 1$ defines a representation   of $\bu$ on $\mathcal{V}_q$. In fact as $\bu$--modules we have
\begin{equation*} V_q(\Lambda_0)\cong \mathcal{V}_q.\end{equation*} 
Further, for all $i\in I$, $u\in\bu^-(0)$, $\eta\in Q$, we have
\begin{equation*} K_i(u\otimes e^\eta) = q^{|\eta|.|\alpha_i|}u\otimes e^\eta,  \ \ \ C(u\otimes e^\eta) = u\otimes e^\eta. \end{equation*}
The highest weight vector in $V_q(\Lambda_0)$ maps to $1\otimes 1$ under
this isomorphism.
\hfill\qedsymbol\end{thm}
Let $\mathcal{V}_\A$ be the image of $V_\A(\Lambda_0)$ under this isomorphism.
Clearly $\mathcal{V}_\A$ is a $\ua$--submodule of $\mathcal{V}_q(\Lambda_0)$
and $\mathcal{V}_q(\Lambda_0)\cong \bq(q)\otimes_{\mathcal{A}}\mathcal{V}_\A$.

Set 
\begin{equation*}\mathcal{L} =\ua^-(0)\otimes \A[Q],\end{equation*}
where $\A[Q]$ is the $\A$--span in $\bq(q)[Q]$  of the elements $e^\eta$.
It follows from  Proposition \ref{heis} that 
\begin{equation*}{\label{lattice}} 
\mathcal{V}_q\cong \bq(q)\otimes_{\mathcal A} \mathcal L.\end{equation*}

We now state our main result. 

\begin{thm}{\label {main}} The lattice $\mathcal{L}$ is preserved by $\ua$,
and 
\begin{equation*} \mathcal{L}\cong \mathcal{V}_\A\end{equation*}
as $\ua$--modules.
\end{thm}
\begin{rem} The case $\frak{g}=sl_2$ was studied in \cite{J3}. In that paper, the author worked over $\mathcal A = \bz[q^{\frac12}, q^{-\frac12}]$ and proved that the corresponding lattice $\mathcal L$ was preserved by $\ua$ and gave the action of the divided powers of the Drinfeld generators on the Schur functions. \end{rem}

The rest of the section is devoted to proving  Theorem \ref{main}.

We begin with the following two lemmas which are  easily deduced from the definition of
$X_i^\pm(z)$ and Proposition \ref{heisr}.
\begin{lem}\label{r=1}
Let $i\in I$, $\eta\in Q$ and $m={|\eta|.|\alpha_i|}$. Then,
\begin{align*}x^+_{i, -m-1}(1\otimes e^\eta) &= \epsilon(\alpha_i,\eta)\otimes  e^{\alpha_i+\eta},\\
x^-_{i, m -1}(1\otimes e^\eta) &= \epsilon(-\alpha_i,\eta)\otimes e^{-\alpha_i+\eta}.
\ \ \ \ \ \ \ \ \hfill\qedsymbol\end{align*}
\end{lem}
\begin{lem}\label{L:product} Let $r,l\in\bz$, $r,l\ge 0$, and let $i,j_1,j_2\cdots j_l\in I$. We have,
\begin{align*} 
&X_i^+(z_1)X_i^+(z_2)\cdots X_i^+(z_r)\left(\tilde\mathcal P^-_{j_1}(w_1)\tilde\mathcal P^-_{j_2}(w_2)\cdots \tilde\mathcal P^-_{j_l}(w_l)\otimes e^\eta\right )\\
&\\
&=\epsilon\cdot\prod_{k=1}^rz_k^{r-k+|\eta|\cdot|\alpha_i|}\prod_{1\le k<s\le r}(f_{i,i}((qz_k)^{-1}, z_s))^{-1}\prod_{1\le k\le r, 1\le s\le l}f_{i,j_s}((qz_k)^{-1}, w_s)\\
&\times \tilde\mathcal P^-_i(z_1)\tilde\mathcal P^-_i(z_2)\cdots \tilde\mathcal P^-_i(z_r)\tilde\mathcal P^-_{j_1}(w_1)\tilde\mathcal P^-_{j_2}(w_2)\cdots \tilde\mathcal P^-_{j_l}(w_l)\otimes e^{r\alpha_i+\eta}\\
&\\
&=\epsilon\cdot (z_1z_2\cdots z_r)^{|\eta|\cdot|\alpha_i|}\prod_{1\le k<s\le r}(z_k-q^{-2}z_s)(z_k-z_s)\prod_{1\le k\le r, 1\le s\le l}f_{i,j_s}((qz_k)^{-1}, w_s)\\
& \times \tilde\mathcal P^-_i(z_1)\tilde\mathcal P^-_i(z_2)\cdots\tilde \mathcal P^-_i(z_r)\tilde\mathcal P^-_{j_1}(w_1)\tilde\mathcal P^-_{j_2}(w_2))\cdots \tilde\mathcal P^-_{j_l}(w_l)\otimes e^{r\alpha_i+\eta},
\end{align*}
where $\epsilon=\ep(r\alpha_i, \eta)\prod_{k=1}^{r-1}\ep(\alpha_i,
k\alpha_i)$. 
\hfill\qedsymbol\end{lem}

Let $\frak{S}_r$ be the symmetric group on $r$ letters and for $\sigma\in\mathfrak S_r$, let $l(\sigma)$ be the length of $\sigma$. 
\begin{lem}\label{id} We have, \begin{equation*}\sum_{\sigma\in\frak S_r}
(-1)^{l(\sigma)}\prod_{k<s}(z_{\sigma(k)}-q^{-2}z_{\sigma(s)}) =
q^{-r(r-1)/2}[r]!\prod_{k<s}(z_k-z_s).\end{equation*} \end{lem}
\begin{proof} Observe that the left--hand side of the equation is an
antisymmetric polynomial in $z_1,z_2,\cdots ,z_r$ and hence is divisible
by the right hand side. Hence, by comparing degrees, we can write,
\begin{equation*}\sum_{\sigma\in\frak S_r}
(-1)^{l(\sigma)}\prod_{k<s}(z_{\sigma(k)}-q^{-2}z_{\sigma(s)}) =
C(q)\prod_{k<s}(z_k-z_s).\end{equation*} But it is easy to see that the
coefficient of $z_1^{r-1}z_2^{r-2}\cdots z_{r-1}$ on the left hand side is
\begin{equation*} \sum_{\sigma\in\frak S_r}q^{-2l(\sigma)}
=q^{-r(r-1)/2}[r]!,\end{equation*} thus proving the proposition.
\end{proof}

\begin{lem}\label {rfact} \begin{enumerate} \item [(i)] Let
$\delta=(\delta_1,\delta_2,\cdots ,\delta_r)\in\bz^r$ be the $r$--tuple
$(r-1,r-2,\cdots ,1,0)$. We have, \begin{equation*}\prod_{j<k}(z_j-z_k)^2
=\left(\sum_\mu a_\mu \sum_{\rho\in\mathfrak
S_r}z_1^{\mu_{\rho(1)}}z_2^{\mu_{\rho(2)}}\cdots
z_r^{\mu_{\rho(r)}}\right),\end{equation*} where the sum is over
$\{\delta+\tau(\delta):\tau\in\mathfrak S_r\}$ and $a_\mu =
(-1)^{l(\tau)}$, if $\mu=\delta+\tau(\delta)$. \item[(ii)]Let $\mathcal R$ be
a commutative ring and let $G\in\mathcal R[[z_1^{\pm1},z_2^{\pm 1},\cdots
,z_r^{\pm1}]]$ be invariant under the action of the symmetric group
$\frak{S}_r$. Then, for all $n\in\bz$, the coefficient of $(z_1z_2\cdots
z_r)^n$ in $\prod_{j<k}(z_j-z_k)^2G$ is divisible by $r!$. \end{enumerate}
\end{lem} \begin{proof} Since, \begin{equation*}\prod_{j<k}(z_j-z_k)
=\sum_{\sigma\in\mathfrak S_r} (-1)^{l(\sigma)}
z_1^{\delta_{\sigma(1)}}z_2^{\delta_{\sigma(2)}} \cdots
z_r^{\delta_{\sigma(r)}},\end{equation*} we get
 \begin{equation*}\prod_{j<k}(z_j-z_k)^2 =\sum_{\sigma,\tau\in\mathfrak
S_r} (-1)^{l(\sigma)+l(\tau)}
z_1^{\delta_{\sigma(1)}+\delta_{\tau(1)}}z_2^{\delta_{\sigma(2)}+\delta_{\tau(2)}}
\cdots z_r^{\delta_{\sigma(r)}+\delta_{\tau(r)}},\end{equation*} which
becomes the formula in (i) on putting $\rho=\sigma\tau$. Part (ii) follows
trivially. \end{proof}

{\it Proof of Theorem \ref{main}}. Using Lemma \ref{L:product} and Lemma \ref{id}, we get
\begin{align*}&\sum_{\sigma\in\frak S_r}X_i^+(z_{\sigma(1)})X_i^+(z_{\sigma(2)})\cdots X_i^+(z_{\sigma(r)}).\left(\mathcal P^-_{j_1}(w_1)\mathcal P^-_{j_2}(w_2)\cdots \mathcal P^-_{j_l}(w_l)\otimes e^\eta\right)\\
&=q^{-r(r-1)/2}[r]!\epsilon\cdot (z_1\cdots z_r)^{|\eta|\cdot|\alpha_i|}\prod_{k<s}(z_k-z_s)^2\prod_{1\le k\le r, 1\le s\le l}f_{i,j_s}((qz_k)^{-1}, w_s)\\
&\times (\tilde\mathcal P^-_i(z_1)\tilde\mathcal P^-_i(z_2)\cdots \tilde\mathcal P^-_i(z_r)\tilde\mathcal P^-_{j_1}(w_1)\tilde\mathcal P^-_{j_2}(w_2)\cdots \tilde\mathcal P^-_{j_l}(w_l))\otimes e^{r\alpha_i+\eta}, \ \ \ \ \ \ (*)\end{align*}
where the constant $\ep$ is defined in Lemma \ref{L:product}.

Set $F=\prod_{k<l}(z_k-z_l)^2$ and let $G$ be the right hand side of (*) divided by $F$. Then  Lemma \ref{rfact} applies, and  
 by collecting the coefficient of $(z_1z_2\cdots z_r)^{-n-1}$ on both 
sides of (*), we find that \begin{equation*}x_{i,n}^{(r)}.
(\tilde\mathcal P^-_{j_1}(w_1)\tilde\mathcal P^-_{j_2}(w_2)\cdots\tilde 
\mathcal P^-_{j_s}(w_l))w_1^{\mu_1}w_2^{\mu_2}\cdots w_l^{\mu_l}
\otimes e^\eta\in \mathcal L,\end{equation*}
for all $\mu_1,\mu_2,\cdots ,\mu_l\in\bz$, $\eta\in Q$, or equivalently that $$(x_{i,n}^+)^{(r)}\mathcal L\subset\mathcal L.$$
One proves similarly that $(x_{i,n}^-)^{(r)}$ preserves $\mathcal L$. 
In particular, by Proposition \ref{adrin} $\mathcal L$ is preserved by $\ua$. 
To complete the proof of the theorem we must prove that 
$$\mathcal L=\mathcal V_\A .$$ 
 Since $1\otimes 1\in\mathcal V_\A$, it follows from  Lemma \ref{r=1} and a simple induction that $e^\eta\in\mathcal V_\A$ for all $\eta\in Q$.
Next, from Theorem \ref{realiz}, we see that  for $i\in I$, $k>0$,
\begin{equation*}\tilde{P}^-_{i,k}(1\otimes e^\eta) = \tilde{P}^-_{i,k}\otimes e^\eta.\end{equation*}
Since by Proposition \ref{heis}, the monomials in the $\tilde{P}^-_{i,k}$'s span $\ua^-(0)$, we see that $\mathcal L\subset \mathcal V_\A$.
The reverse inclusion $\mathcal V_\A\subset \mathcal L$ is now clear, for 
$$\mathcal V_\A =\tilde\ua(1\otimes 1)\subset\mathcal L $$
since $1\otimes 1\in\mathcal L$. \hfill\qedsymbol

\section{Specialization to a root of unity} Throughout this section, we let $N$ denote the Coxeter number of $\frak{g}$. It is well--known \cite{Bo} that 
$N= n+1$  (resp. $2n-2$, 12, 18, 30) 
if $\frak g$ is of type $A_n$ (resp. $D_n$, $E_6$, $E_7$, $E_8$). 
Let $ \zeta\in\bc^*$  denote  a primitive $l^{th}$ root of unity, where $l$ is a 
non--negative integer coprime to $N$. 
Set $n=|I|$. Finally, for any 
$g\in\A$ we let 
$g_\zeta\in\bc^*$ be the element obtained by setting $q=\zeta$. 
\begin{lem}\label{det} Let $[A]$ denote the $n\times n$--matrix with coefficients in $\A$ whose $(i,j)$--th entry is $[a_{ij}]$. Then,
\begin{align*} {\text{det}}[A] &= [n+1],\ \ \ {\text{if $\frak g$ is of type $A_n$}},\\
&= [2](q^{n-1}+q^{n-1}), \ \ \ {\text{if $\frak g$ is of type $D_n$}},\\
&=(q^4+q^{-4}-1)(q^2+q^{-2}+1), \ \ \ {\text{if $\frak g$ is of type $E_6$}},\\
&=[2](q^6+q^{-6}-1), \ \ \ {\text{if $\frak g$ is of type $E_7$}},\\
&= q^8+q^6+q^{-6}+q^{-8}-q^2-1-q^{-2} \ \ \ {\text{if $\frak g$ is of type $E_8$}}.\end{align*}
Further for all $k>0$, we have  \begin{equation*}({\text{det}}[A])_{\zeta^k} = {\text{det}}[A]_{\zeta^k} \ne 0.\end{equation*}
\end{lem}
\begin{pf} The calculation of the determinant is straightforward.
 If $\frak g$ is of type $A_n$, then it is easy to see that  for all $k>0$, $${\text{either}}\ \ \zeta^{2k} =1,\ \ \text{or}\ \  \zeta^{2k(n+1)} \ne 1.$$  This proves the second statement of the Lemma for $\frak g$ of type $A_n$. The other cases are 
proved by a similar analysis: in the hardest case  $E_8$ one checks that $q^8+q^6+q^{-6}+q^{-8}-q^2-1-q^{-2}$ divides $q^{60}-1$ in $\A$. The result follows.
\end{pf}

Let $\bc_{\zeta}$ be the
one--dimensional $\mathcal A$--module defined by sending
$q\to\zeta$.  Let $\bu_\zeta$ be the algebra over $\bc$ defined by,
\begin{equation*}\bu_\zeta=\ua\otimes_\A \bc_\zeta.\end{equation*}
The subalgebras $\bu^\pm_\zeta$ and $\bu_\zeta^\pm(0)$ of $\bu_\zeta$ are defined in the obvious way and we have
$$\bu_\zeta=\bu_\zeta^-\bu_\zeta^0\bu_\zeta^+.$$ Given an element $u\in\ua$, we denote by $u$ the element $u\otimes 1$ in $\bu_\zeta$. It follows from Proposition  \ref{heisr}  that we have a representation $\pi_\zeta: \bu_\zeta(0)\to  
{\text{End}}(\bu_\zeta^-(0))$.

\begin{prop} {\label{heiszeta}} 
\begin{enumerate}
\item[(i)] For all $i\in I$, $k\in\bz$, $k>0$ there exist elements $h^{i,k}\in\bu_\zeta(0)$ such that
\begin{align*} [h^{i, k}, \tilde h_{j,m}] &= \delta_{k,-m}\delta_{i,j},\\
 [h^{i,k}, \tilde P^-_{j,m}]&=\delta_{k,m}\delta_{i,j}.\end{align*}
\item[(ii)] $\pi_\zeta$ is an irreducible representation of $\bu_\zeta(0)$. \end{enumerate}
\end{prop}
\begin{proof}  For $k\in\bz$, $k>0$, we know by Lemma \ref{det} that  the matrix  $[A]_{\zeta^k}$ is invertible. Let $b_{ij}(k)$ denote the inverse of this matrix. 

For $i\in I$, $k\in\bz$, $k>0$, set $$h^{i,k} =\sum_{j\in I}b_{ij}(k) \tilde h_{j,k}. $$
Clearly $h^{i, k}$ satisfies \begin{equation*} [h^{i, k}, \tilde h_{j,m}] = \delta_{k,-m}\delta_{i,j}.\end{equation*}
The second formula in (i) is now clear from Proposition \ref{heis}.

To prove (ii), assume that $W$ is a  submodule of $\bu^-_\zeta(0)$ and let $0\ne w\in W$. By Proposition \ref{heis}, we can choose $i\in I$, $k\in \bz$, $k>0$, such that  
$$w=\sum_{ r=0}^n (\tilde P^-_{i,k} )^r w_r$$
where $w_r$ is a polynomial in the elements $\tilde P^-_{j,l}$, $j\ne i$, $1\le l\le k$ and $\tilde P^-_{i,l}$, $1\le l<k$. Applying $h^{i,k}$ to $w$ repeatedly we find that 
$w_n\in W$. Repeating the argument we find that $1\in W$ thus proving the Proposition.
\end{proof}

We now turn to the representations of $\bu_\zeta$. Given $\Lambda=\sum_{i\in\hat I}n_i\Lambda_i$, $n_i\ge 0$, set, 
$$W_\zeta(\Lambda) = V_\A(\Lambda)\otimes_\A\bc_\zeta .$$
It follows from Proposition \ref{vlatt}  that 
$W_\zeta(\Lambda)$ is a representation of $\bu_\zeta$. Again, for $v\in V_\A(\Lambda)$, we let $v\in W_\zeta(\Lambda)$ be the element $v\otimes 1$. 
Clearly $\bu_\zeta^+.v =0$ and  
\begin{equation*} W_\zeta(\Lambda) = \bu_\zeta. v_\Lambda.\end{equation*}
Set $W_\zeta(\Lambda)_{\mu,n} =(V_\A(\Lambda)\cap V_q(\Lambda)_{\mu,n})\otimes_\A \bc_\zeta$. 
Then one knows from \cite{L, L2} that, 
\begin{equation*} W_\zeta(\Lambda)=\bigoplus_{\mu,n} W_\zeta(\Lambda)_{\mu, n}, \ \ {\text{dim}}_\bc W_\zeta(\Lambda) ={\text{dim}}_{\bq(q)} V_q(\Lambda), \end{equation*}
and $w\in W_\zeta(\Lambda)_{\mu,n}$ iff
\begin{align*} K_i.v = \zeta^{\mu_i'}v, \ \ & \ \ \genfrac{[}{]}{0pt}{}{K_i,0}{l}.v = \mu_i''v,\\
D.v =\zeta^{n'}v,\ \ &\ \ \genfrac{[}{]}{0pt}{}{D,0}{l}.v = n''v,\end{align*}
where $\mu_i=\mu_i'+l\mu_i''$, $0\le \mu_i'<l$, and $n'$ and $n''$ are defined similarly. 

Turning now to the level one basic representation, we see from Theorem \ref{main} that \begin{equation*}
W_\zeta(\Lambda_0)\cong \bu_\zeta^-(0)\otimes \bc_\zeta[Q].\end{equation*}
The main result of this section is:
\begin{thm} 
 $W_\zeta(\Lambda_0)$ is irreducible.
\end{thm}

\begin{pf} Set $W=W_\zeta(\Lambda_0)$ and  let $0\ne W'$ be a submodule of $W$.
Then $W'$ contains a non--zero vector $w\in W_{\mu,n}$ such that 
$$\bu_\zeta^+.w =0.$$ It is clear from Theorem \ref{realiz} that $w$ must be of the form $w_\mu\otimes e^\mu$ for some $w_\mu\in\bu_\zeta^-(0)$ with $$\bu^+_\zeta(0).w_\mu =0.$$ By Proposition 
\ref{heiszeta} we see that this forces $w_\mu =1$ and hence that $1\otimes e^\mu\in W'$. Proposition \ref{r=1} now shows that $1\otimes e^{\nu }$ for all $\nu\in Q$ and hence finally that $W'=W$. \end{pf}

\end{document}